\documentclass{amsart}
\usepackage{amsfonts,amssymb,latexsym,amsmath,amsthm}
\usepackage{url}
\usepackage{enumerate}
\usepackage{graphicx}
\usepackage{epstopdf}
\def\la{{\lambda}}

\def\la{{\lambda}}
\def\Bop{{\mathbb B}}
\def\Cop{{\mathbb C}}
\def\Hop{{\mathbb H}}
\def\Sop{{\mathbb S}}
\def\ZZ{{\mathbb Z}}
\def\Ht{{\tilde H}}
\def\coeff{{\Big|}}

\def\doff{{\rm{doff}}}
\def\dinv{{\rm{dinv}}}
\def\area{{\rm{area}}}
\def\Des{{\text{Des}}}
\def\touch{{\rm{touch}}}

\def\WP{{\mathcal W}P}

\newcommand{\rev}[1]{\overleftarrow{#1}}
\newcommand{\pchoose}[2]{\begin{pmatrix}#1\\ #2\end{pmatrix}}
\newcommand{\bchoose}[2]{\begin{bmatrix}#1\\ #2\end{bmatrix}}

\newtheorem{thrm}{Theorem}[section]
\newtheorem{lem}[thrm]{Lemma}
\newtheorem{prop}[thrm]{Proposition}
\newtheorem{cor}[thrm]{Corollary}
\newtheorem{conj}[thrm]{Conjecture}
\theoremstyle{definition}

\newtheorem{remark}[thrm]{Remark}
\newtheorem{example}[thrm]{Remark}
\numberwithin{equation}{section}

\newdimen\squaresize \squaresize=12pt
\newdimen\thickness \thickness=0.4pt

\def\square#1{\hbox{\vrule width \thickness
     \vbox to \squaresize{\hrule height \thickness\vss
        \hbox to \squaresize{\hss#1\hss}
     \vss\hrule height\thickness}
\unskip\vrule width \thickness}
\kern-\thickness}

\def\vsquare#1{\vbox{\square{$#1$}}\kern-\thickness}

\def\thisbox#1{\kern-.09ex\fbox{#1}}
\def\downbox#1{\lower1.200em\hbox{#1}}
\title
{
Dyck Paths with Forced and Forbidden Touch Points and
$q,t$-Catalan building blocks}
\author{J. Haglund}
\address{Department of Mathematics, University of Pennsylvania, Philadelphia, PA 19104-6395}
\email{jhaglund@math.upenn.edu }

\author{J. Morse}
\address{Department of Mathematics, Drexel University, Philadelphia, PA 19104}
\email{morsej@math.drexel.edu}

\author{M. Zabrocki}
\address{Department of Mathematics and Statistics, York University, Toronto, Ontario M3J 1P3, CA}
\email{zabrocki@mathstat.yorku.ca}

\thanks{Work supported by NSF grants DMS-0553619 and DMS-0901467}


\begin{document}
%\numberwithin{theorem}{section}

%\title[Dyck Paths with Forced and Forbidden Touch]
%{
%%Hall-Littlewood Symmetric Functions indexed by Compositions and
%Dyck Paths with Forced and Forbidden Touch Points and
%$q,t$-Catalan building blocks}
\title[A compositional shuffle conjecture]{A compositional Shuffle conjecture
specifying touch points of the Dyck path}

\begin{abstract}
We introduce a $q,t$-enumeration of Dyck paths which are forced to touch the main diagonal 
at specific points and forbidden to touch elsewhere 
and conjecture that it describes the action of
the Macdonald theory $\nabla$ operator applied to a Hall-Littlewood 
polynomial.  Our conjecture refines several earlier conjectures concerning
the space of diagonal harmonics including the ``shuffle conjecture"
(Duke J. Math. $\mathbf {126}$ ($2005$), pp. $195-232$) for $\nabla e_n[X]$.
We bring to light that certain generalized Hall-Littlewood polynomials
indexed by compositions are the building blocks for the algebraic
combinatorial theory of $q,t$-Catalan sequences and we prove a number of 
identities involving these functions.
%and connect them to decompositions
%of $e_n(X)$ introduced by Garsia and Haglund (Discrete Math. $\mathbf 256$ ($2002$),
%pp. $677-717$) in their study of the $q,t$-Catalan sequence. 
\end{abstract}
\maketitle

\section{Introduction}
Our study concerns the combinatorics behind the character of the space 
of diagonal harmonics $\text{DH}_n$ and identities involving Macdonald 
polynomials that can be used to form expressions for this character.  
%Background references for these topics
%are the papers \cite{Hai94},\cite{GaHm96b},\cite{GaHa01},\cite{GaHa02} and 
%\cite{HHLRU05}.  The results in these papers and others are summarized in the book
%\cite{Hag08}.  A brief synopsis of the current state of affairs is as follows.
At the root of this theory is a linear operator $\nabla$, 
introduced in \cite{BGHT99}, under which the
modified Macdonald polynomials $\tilde H _{\mu}[X;q,t]$ are 
eigenfunctions.  
%*%Haiman 
\cite{Hai02} proved that the Frobenius image of the 
character of $\text{DH}_n$ equals
$\nabla e_n$.  This gives an explicit expression involving
rational functions in $q,t$ for the multiplicity of an irreducible
indexed by a partition $\lambda$ in the character of $\text{DH}_n$.  

An important open problem in this area is the ``shuffle conjecture" 
of \cite{HHLRU05} which asserts that the coefficient of $m_{\lambda}$ 
in $\nabla e_n$ simplifies to 
a $q,t$ statistic on lattice paths.  A major breakthrough in this direction 
was made with the conjectured combinatorial formula of \cite{Hag03} for the 
coefficient of $m_{1^n}=s_{1^n}$ in $\nabla e_n$.  In this case, the 
coefficient is a bi-graded version of the sign character and it is
called the $q,t$-Catalan $C_n(q,t)$ since it reduces to the $n$th Catalan 
number when $q=t=1$.  The combinatorial formula for $C_n(q,t)$ was
proven in \cite{GaHa01},\cite{GaHa02} and pursuant work \cite{Hag04a}
also settled the shuffle conjecture for partitions of hook-shape.
However, the general case remains a mystery.

An unrelated study of Macdonald polynomials led to the discovery
\cite{LLM03}
of a new family of symmetric functions called $k$-Schur functions 
$s_\lambda^{(k)}[X;t]$ which were conjectured 
to refine the special combinatorial properties held by Schur functions.
The $k$-Schur functions have a number of conjecturally equivalent 
characterizations and it has now been established in \cite{LaMo08a,Lam06b} 
that those introduced in \cite{LM07}
refine combinatorial, geometric and representation 
theoretic aspects of Schur theory.   This prompted Bergeron, Descouens, 
and Zabrocki to explore the role of $k$-Schur functions in the 
$q,t$-Catalan theory.  To this end, they conjectured in \cite{BDZ10} 
that the coefficient of $s_{1^n}$ in $\nabla s_{1^n}^{(k)}[X;t]$ is 
a positive polynomial in $q,t$ and they proved their conjecture
for the case that $t=1$.

Our work here was initially motivated by a desire to find a combinatorial
description for this coefficient in general, ideally in terms of a 
$q,t$-statistic on lattice paths as with the $q,t$-Catalans.
%in \cite{LLM03} and studied further in \cite{LaMo03a},\cite{LaMo03b},\cite{LaMo08a} 
%and in a number of
%more recent papers by a variety of authors.  
%%They reduce to Schur functions when $k\ge n$.)  
We found such a description, but more remarkably this led us to
discover that a natural setting for the combinatorial theory
of $DH_n$ is created by applying $\nabla$ to the general
set of Hall-Littlewood polynomials indexed by compositions.
To be precise, it was proven in \cite{LaMo03a} that the
$k$-Schur function $s_{1^n}^{(k)}[X;t]$ is merely a 
certain Hall-Littlewood polynomial.  This led us to 
study $\nabla$ on a Hall-Littlewood polynomial indexed by 
any partition $\lambda$.  But in fact, our work carries 
through to the family of polynomials $C_\alpha[X;q]$, for any composition
$\alpha$, defined in terms of operators similar to Jing operators.

A key component in the proof of the $q,t$-Catalan 
conjecture \cite{GaHa02} is the use of symmetric functions $E_{n,k}[X;q]$ 
which decompose $e_n$ into pieces that remain positive under 
the action of $\nabla$.  We have discovered that the $C_\alpha[X;q]$ can be used as
building blocks in the $q,t$-Catalan theory that 
decompose the $E_{n,k}[X;q]$ 
into finer pieces, still positive under the action of $\nabla$. 
Our conjectures on these building blocks thus 
refine earlier conjectures
involving $E_{n,k}[X;q]$,  
the conjectures in \cite{BDZ10}, 
the shuffle conjecture, and the conjectures in \cite{BGHT99} 
asserting that $\nabla$ applied to Hall-Littlewood functions have $q,t$-positive 
Schur coefficients.  
Loehr and Warrington \cite{LoWa08} introduced an 
intricate conjecture for the combinatorics of $\nabla$ applied 
to a Schur function $s_{\lambda}$.  Our conjecture is extremely simple,
describes the action of $\nabla$ on a larger set of symmetric functions
than just a basis, and refines the conjecture of Loehr and Warrington when
$\nabla$ acts on the Schur function $s_{(n-k,1^k)}$ \cite[Conjecture 3]{LoWa08} 
as explained at the end of Section \ref{Nab}.

Garsia, Xin, and Zabrocki \cite{GXZ10} using work of Hicks \cite{Hic10} have now proven 
our generalized $q,t$-Catalan conjecture and expanded
the result giving a ``compositional $q,t$-Schr\"oder''
theorem.

\section{Definitions and notation}

\subsection{Combinatorics}

A Dyck path is a lattice path in the first quadrant of 
the $xy$-plane from the point $(0,0)$ to
the point $(n,n)$ with steps $+(0,1)$ and $+(1,0)$
which stays above the line $x=y$.   
For a Dyck path $D$, the cells in the $i^{th}$ row are
those unit squares in the $xy$-plane that are below the path 
and fully above the line $x=y$ whose NE corner has a $y$ coordinate
 of $i$.  The set of Dyck paths from $(0,0)$
 to $(n,n)$ will be denoted $DP^n$ and the number of
 paths in this set is well known to be the Catalan number
 $$C_n = \frac{1}{n+1} \binom{2n}{n}~.$$

For a Dyck path $D$, let $a_i=a_i(D)$ equal the number of cells in the $i^{th}$ 
row of $D$.  It is always true $a_1 = 0$ and $0 \leq a_{i+1} \leq a_i + 1$.
We define the arm sequence $\text{arm}(D) = (a_1, a_2, \ldots, a_n)$
and note this completely determines $D$.  We consider two statistics 
(non-negative integers) on Dyck paths.  The $area$ statistic is the number 
of whole cells which are below the path and above the diagonal,
or
$$
\area(D) = \sum_{i=1}^n a_i.
$$ 
The $dinv$ statistic
%, due to Mark Haiman,
is defined as 
$$
\dinv(D) = \sum_{1 \leq i < j \leq n} \chi( a_i - a_j \in \{ 0, 1\})
$$
where $\chi( \text{true} ) =1$ and $\chi(\text{false}) = 0$.  

\begin{example}
The Dyck path $D$ with arm sequence $(0^n)$ has 
$\area(D) = 0$ and $\dinv(D) = \binom{n}{2}$.  The Dyck path $D'$ with
arm sequence $(0, 1, 2, 3, \ldots, n-1)$ has $\area(D') = \binom{n}{2}$
and $\dinv(D') = 0$.

\begin{center}
\begin{tabular}{ccccc}&
\includegraphics[width=.75in]{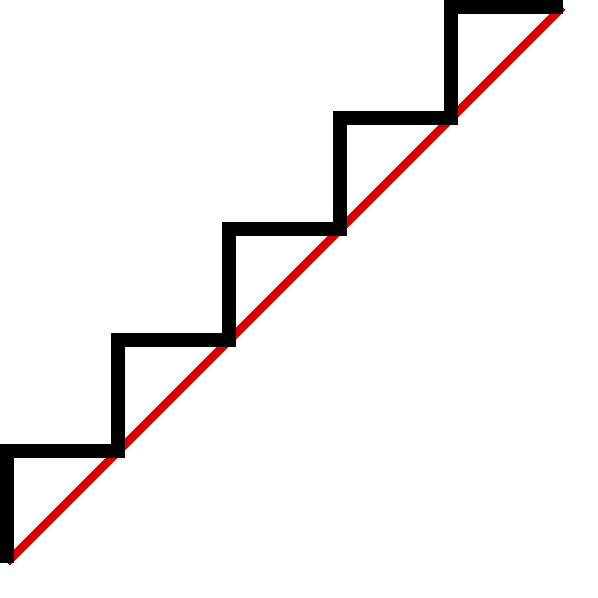}&
\includegraphics[width=.75in]{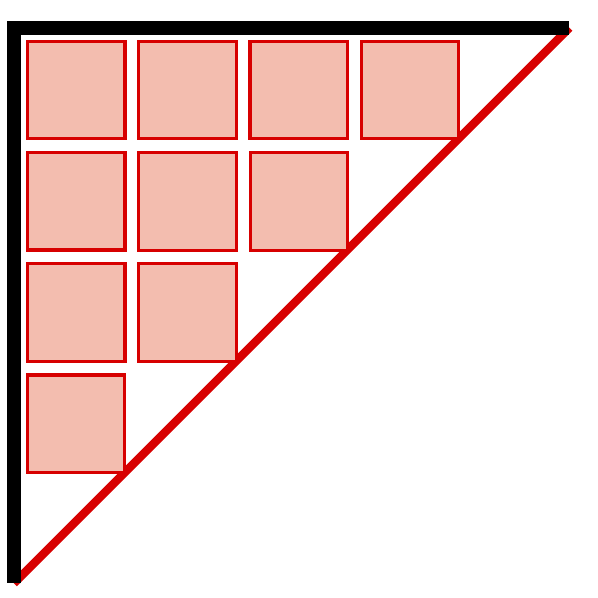}&
\includegraphics[width=.75in]{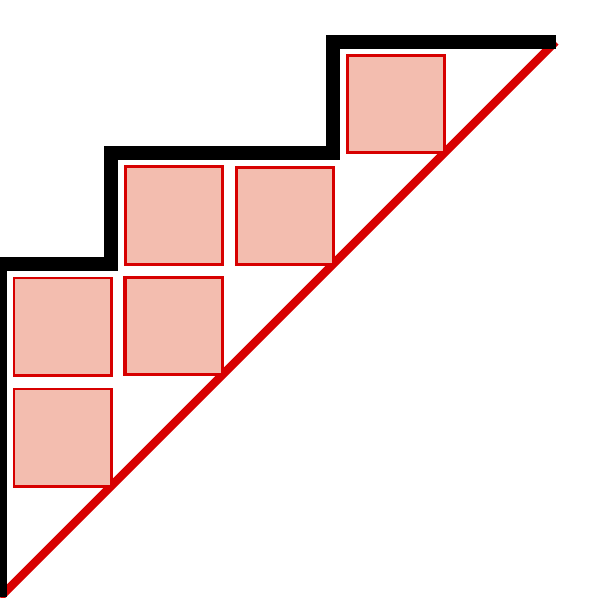}&
\includegraphics[width=.75in]{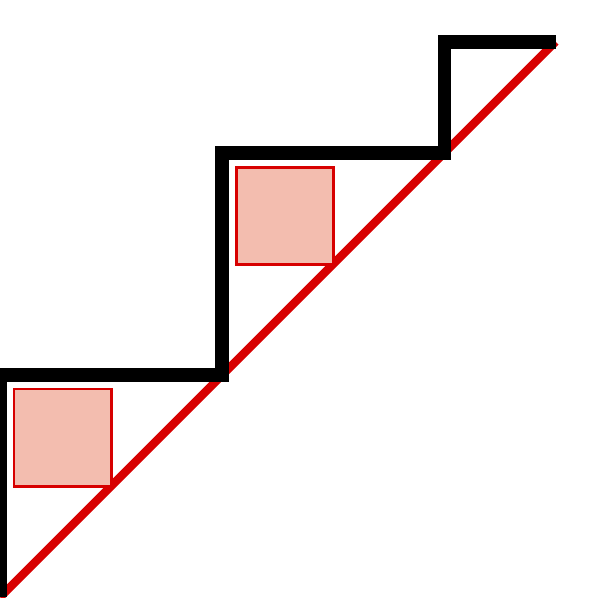}\\
arm sequence &(0,0,0,0,0)&(0,1,2,3,4)&(0,1,2,2,1)&(0,1,0,1,0)\\
\area&0&10&6&2\\
\dinv&10&0&4&7
\end{tabular}
\end{center}
\end{example}

\begin{remark} 
The original proof of the combinatorial interpretation
of the $q,t$-Catalan was stated in terms of a third statistic $bounce(D)$.
Since we are more cleanly able to formulate our results in terms of the $\dinv(D)$
statistic we choose to state all the results in this paper in
terms of the $\dinv(D)$ statistic however the reference \cite[p. 50]{Hag08}
describes an automorphism $\phi$ on $DP^n$ such that 
$\area(\phi(D)) = \text{bounce}(D)$
and $\dinv(\phi(D)) = \area(D)$.
\end{remark}

%\todo{Should figure out what the \touch composition corresponds to
%with the map that exchanges \dinv -> \area, \area -> bounce,
%I kind of did this and it doesn't seem like a useful statistic.  One
%can figure out what the sizes of the parts of the \touch composition
%by looking at lengths of certain horizontal (or vertical steps are).
%}

We make use of a partial order on Dyck paths; namely $D_1 \leq D_2$ if 
$\text{arm}(D_1) \leq \text{arm}(D_2)$, component-wise.
In this case we say that $D_1$ is `below'
$D_2$ because $D_1$ will not cross $D_2$ and is hence weakly `between'
$D_2$ and the diagonal.

A composition $\alpha$ of $n$, denoted $\alpha\models n$, is an
integer sequence $\alpha = (\alpha_1, \alpha_2, \ldots, \alpha_r)$ 
with $\alpha_i \geq 1$ and where 
$|\alpha| = \alpha_1+ \alpha_2+ \cdots + \alpha_r = n$.   
The length of $\alpha$ is $\ell(\alpha)=r$.
We shall also use ${\rev{\alpha}} = (\alpha_{\ell(\alpha)},
\alpha_{\ell(\alpha)-1}, \ldots, \alpha_2, \alpha_1)$.  
For any composition $\alpha$, we define
\begin{equation}
n(\alpha) = \sum_{i=1}^{\ell(\alpha)} (i-1) \alpha_i\,.
\end{equation}
The descent set of a composition $\alpha$ is defined to be
$$\Des(\alpha) = 
\{\alpha_1 , \alpha_1 + \alpha_2 , \ldots, \alpha_1 + \alpha_2 + \cdots + \alpha_{\ell(\alpha)-1} \}~.
$$
There is a common partial order defined on composition $\alpha, \beta \models n$
by letting $\alpha \leq \beta$ 
when $\alpha$ is `finer' than $\beta$, i.e. $\Des(\beta) 
\subseteq \Des(\alpha)$. 
If $\alpha$ is a composition of $n$, $DP(\alpha)$ 
represents the Dyck path consisting of $\alpha_1$ 
steps in the North $(0, 1)$ direction followed by $\alpha_1$ steps in the 
East $(1, 0)$ direction, $\alpha_2$ $(0,1)$ steps 
followed by $\alpha_2$ $(1,0)$ steps, etc.

A partition $\lambda = (\lambda_1,\ldots,\lambda_r)$ is a non-increasing 
sequence of positive integers.
When $\lambda$ is a partition of $n$, denoted $\lambda\vdash n$,
$|\lambda|=\sum \lambda_i = n$.  The length of $\lambda$ is
$\ell(\lambda)=r$.
Given a partition $\lambda$, we set
\begin{equation}
m(\lambda)= (m_1(\lambda), m_2(\lambda), m_3(\lambda), \cdots,
m_{|\la|}(\lambda))\,,
\end{equation}
where the numbers $m_i(\lambda)$ represent the number of parts of size $i$
in $\lambda$.  The conjugate of a partition $\lambda$ is the
partition $\lambda'= (\lambda_1',\lambda_2',\ldots,\lambda_m')$ where
$\lambda_i'$ is the number of parts of $\lambda$ that are
at least $i$.  
Partitions are generally considered
to be compositions with parts arranged in non-increasing order.
Hence, notions defined on compositions apply to partitions as well.
Generally, we will use the symbols $\alpha, \beta, \gamma$ to represent 
compositions and $\la, \mu, \nu$ for partitions.  

For a given Dyck path $D$, $\touch(D)$ denotes  the composition 
$\gamma = (\gamma_1, \gamma_2, \ldots, \gamma_{\ell(\gamma)}) \models n$ 
that specifies which rows the Dyck path `touches' the diagonal.  
That is, for $\text{arm}(D) = (a_1, a_2, \ldots, a_n)$, $a_k = 0$ if and only if
$k=1$ or $k-1 \in \Des(\gamma)$.
The title of this paper comes from the notion of the $\touch$ composition. By 
requiring that $\touch(D) = \alpha$ for a fixed composition $\alpha$,
we have specified that the Dyck path will touch the diagonal in rows
$1, 1+\alpha_1, 1+\alpha_1+\alpha_2, \ldots, 1+\alpha_1+\alpha_2+ \cdots+ 
\alpha_{\ell(\alpha)-1}$ and is forbidden to touch the diagonal in the other
rows.  Note that under this definition, we view all paths as touching the 
diagonal in row $1$, but none touching in row $n+1$, and we say the path 
touches the diagonal $\ell (\alpha)$ times.
The partial order on compositions is consistent with the partial order on
Dyck paths in the sense that if $D_1$ and $D_2$ are Dyck paths such that
$D_1 \leq D_2$, then $\touch(D_1) \leq \touch(D_2)$.

Using these notions we introduce a new statistic $\doff_\alpha(D)$ 
for a given Dyck path $D$ with $\touch(D)\leq\alpha$. 
% This statistic represents a {\it \dinv offset}
%in some of the equations where we conjecture or provide
%a combinatorial interpretation. 
If $\text{arm}(D)=(a_1, a_2, \ldots, a_n)$,
let $r_1$ be the number of rows such that $a_i=0$ for 
$1\leq i \leq \alpha_1$, $r_2$ be the number of rows such that $a_i=0$ for 
$\alpha_1< i \leq \alpha_1+\alpha_2$, and more generally
$r_k = \#\{ i : a_i = 0\hbox{ and }\sum_{j=1}^{k-1} \alpha_j < i \leq\sum_{j=1}^{k} \alpha_j\}$.
We then set
\begin{equation}\label{doffstat}
\doff_\alpha(D) = \sum_{k=1}^{\ell(\alpha)}(\ell(\alpha)-k)r_k~.
\end{equation}

\begin{example}
If $\text{arm}(D) = (0,1,2,0,1,2,2,1,0,1,2,3,2,1)$, then
$\touch(D) = (3,5,6)$.
Taking $\alpha = (8,6)$ we have $\doff_{(8,6)}(D) = 2$.

\begin{center}
\includegraphics[width=2in]{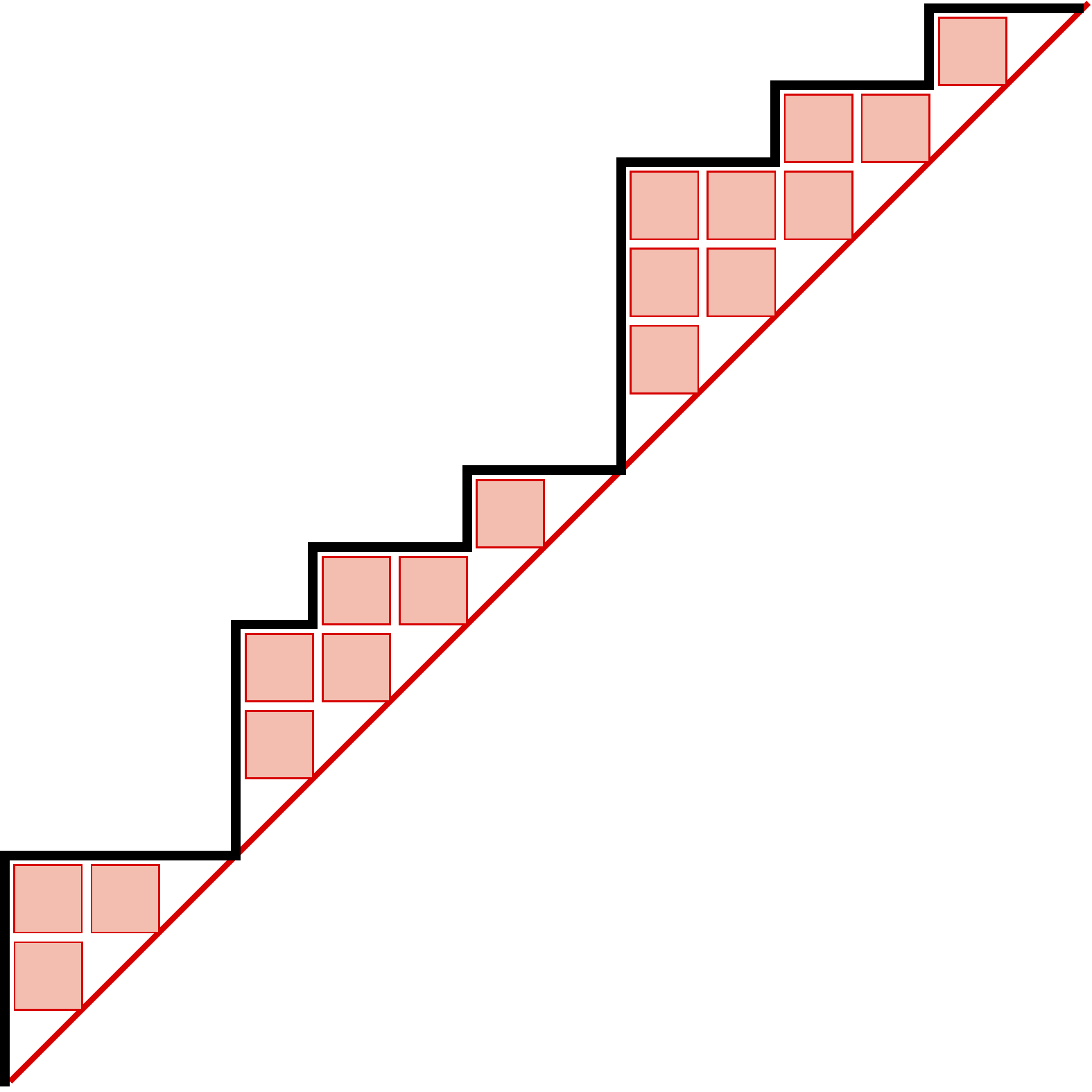}
\end{center}

\noindent
\end{example}

The only Dyck path with $\touch(D) = (1^n)$ has $\text{arm}(D) = (0^n)$.  
There are $C_{n-1}$ Dyck paths with $\touch(D) = (n)$ and more generally 
there are $\prod_{i=1}^{\ell(\alpha)} C_{\alpha_i-1}$ Dyck paths such 
that $\touch(D) = \alpha$.  Note that if $D$ and $E$ are two Dyck paths 
with $\text{touch}(D)=\text{touch}(E)\le \alpha$, then
$\text{doff} _{\alpha} (D)=\text{doff} _{\alpha} (E)$.

The results we mention so far are stated in terms of Dyck paths, but 
we will require the notion of parking functions to state the generalization 
of the shuffle conjecture. For a Dyck path $D$ in $DP^n$ with 
$\text{arm}(D)=(a_1,\ldots,a_n)$, let $\WP_D$ be the set of words of 
length $n$ in the alphabet $\{1, 2, \ldots , n\}$ such that $w_i < w_{i+1}$ 
if $a_i < a_{i+1}$.
We use the notation $x^w$ to denote the monomial 
$x_{w_1} x_{w_2} \cdots x_{w_n}$. We also 
define an extension of the $\dinv$ statistic for words in 
$\WP_D$ by setting 
\begin{align*}
\dinv(w) = |\{&(i,j):1 \leq i < j \leq n, a_i = a_j\hbox{ and } w_i < w_j \}|\\
&+|\{(i,j): 1 \leq i < j \leq n, a_i = a_j + 1\hbox{ and }w_i > w_j\}|~.
\end{align*}

\subsection{Symmetric functions}
Let $X$ represent a sum of an infinite set of
variables $X = x_1 + x_2 + x_3 + \cdots$ considered as
elements of the ring of polynomial series in an
infinite number of variables of bounded degree.  
For $r >0$, let $p_r$ represent a linear and algebraic 
morphism which acts on polynomial series by $p_r[x] = x^r$.
That is for two polynomial series of bounded degree $A$ and $B$,
\begin{align*}
p_r[A+B] &= p_r[A] + p_r[B]\\
p_r[A - B] &= p_r[A] - p_r[B]\\
p_r[AB] &= p_r[A] p_r[B]\\
\end{align*}  and in particular,
$p_r[X] = x_1^r + x_2^r + x_3^r + \cdots$ represents the
$r^{th}$ power sum in the 
variables $\{x_1, x_2, x_3, \ldots\}$. 
The ring of symmetric functions over the field $F$
is defined to be the polynomial
ring
$$\Lambda = F[ p_1[X], p_2[X], p_3[X], \ldots]~.$$
For our purposes, we choose the field $F$ to be the ring of rational power 
series in the variables $q,t,u,z$ over $\mathbb{Q}$ where each of the
parameters $q,t,u,z$ all have the property that
$p_r[a] = a^r$ for each $a = q,t,z,u$.

Generally our symmetric functions $f$ will be considered as
polynomials in the elements $p_r$ and then the notation $f[A]$
represents $f$ with each $p_r$ replaced by $p_r[A]$.  The degree
of $p_r$ is $r$ and the degree of a symmetric function $f$ is determined
by the degree of the monomials in the power sums which appear in $f$.
Following the notation of Macdonald \cite{Macdonald}, we have the 
power sum basis $p_\lambda[X]$, Schur basis $s_\lambda[X]$,
homogeneous basis $h_\lambda[X]$, and elementary basis $e_\lambda[X]$.

In the expressions of variables it is useful to have
a special symbol $\epsilon$ which will represent
a value of negative one but behaves differently than a negative symbol. 
If $f$ is of homogeneous degree $r$,
\begin{align*}
f[\epsilon X] &= (-1)^r f[X]\\
f[-\epsilon X]&= \omega (f[X])\\
\end{align*}
where $\omega$ is an involution on symmetric functions such that
$\omega( p_\lambda[X] ) = (-1)^{|\la|+\ell(\la)}p_\lambda[X]$,
$\omega( e_n[X] ) = h_n[X]$ and $\omega( s_\lambda[X] ) = s_{\lambda'}[X]$.
We will also make use of the standard Hall scalar product
which is defined by
$$\left< p_\la[X], p_\mu[X]/z_\mu \right> = \left< s_\la[X], s_\mu[X] \right> = \chi(\la = \mu)$$
where $z_\mu = \prod_{i \geq 1} m_i(\mu)! i^{m_i(\mu)}$.

For any symmetric function $f$, multiplication by $f$ is an operation on symmetric functions
which raises the degree of the symmetric function by $deg(f)$.  If we define
$f^\perp$ to be the operation which is dual to multiplication in the sense that
$$\left< f^\perp g, h \right> = \left< g, f h \right>,$$
then $f^\perp$ is an operator which lowers the degree of the symmetric function
by $deg(f)$.  It is not difficult to show that
\begin{equation}
\label{hkperp}
f[X + z] = \sum_{k\geq0} z^k (h_k^\perp f)[X],
\end{equation}
\begin{equation}
\label{ekperp}
f[X - z] = \sum_{k\geq0} (-z)^k (e_k^\perp f)[X].
\end{equation}

In addition we will refer to the form of the
Macdonald basis $\Ht_\lambda[X;q,t]$ that
is relevant to the study of the $n!$ Theorem \cite{Hai01} and the
$q,t$-Catalan numbers.  The relation of this basis to
the integral form $J_\mu[X;q,t]$ of \cite{Macdonald} is
$$\Ht_\mu[X;q,t] = t^{n(\mu)} J_\mu\left[ \frac{X}{1-1/t};q,1/t \right].$$
It is also characterized as the unique basis such that
$$\left< \Ht_\mu[X(1-1/t);q,t], \Ht_\la[X(1-q);q,t] \right> = 0$$ if
$\la \neq \mu$ and $\left< \Ht_\mu[X;q,t], h_n[X] \right> = 1$.

We are particularly interested in the Hall-Littlewood symmetric functions.
Following the notation of Macdonald we define the functions
$Q'_\la[X;q]$ to be the basis of the symmetric functions which satisfy
$$\left< Q'_\lambda[X(1-q);q], Q'_\mu[X;q] \right> = 0$$
if $\la \neq \mu$ and $\left< Q'_\la[X;q], h_n[X] \right>  = q^{n(\lambda)}$.
Relating the definitions of the Hall-Littlewood and Macdonald symmetric
functions, we note that
$$Q'_\la[X;q] = \Ht_\la[X;0,1/q] q^{n(\la)}= \sum_{\la} K_{\la\mu}(q) s_\la[X]~.$$

The operator $\nabla$ was introduced in \cite{BGHT99}
and is defined by
\begin{equation}\label{nabladef}
\nabla \Ht_\la[X;q,t]  = t^{n(\lambda)} q^{n(\la')} \Ht_\la[X;q,t]~.
\end{equation}
This operator has been fundamental to the study of the $q,t$-combinatorial
identities associated with $\text{DH}_n$ and Macdonald polynomials.  Its 
definition is chosen so that
\begin{equation}\label{defqtcat}
\left< \nabla( e_n[X] ), e_n[X] \right> = C_n(q,t)
\end{equation}
where $C_n(q,t)$ is the $q,t$-Catalan polynomial.  
%*%Garsia and the first author 
References \cite{GaHa01, Hag03, Hag08} showed that
\begin{equation}\label{CIqtcatalan}
C_n(q,t) = \sum_{D \in DP^n} t^{\area(D)} q^{\dinv(D)}
\end{equation} with the sum
over all Dyck paths of length $n$.

\begin{example}
A small example is $C_3(q,t) = q^3 + qt+ q^2t + qt^2 + t^3$,
whose terms can be computed (in order) from the 5 Dyck paths of length $3$ 
with respective arm sequences $(0,0,0)$, $(0,0,1)$, $(0,1,0)$, $(0,1,1)$ 
and $(0,1,2)$:  
\begin{center}
\includegraphics[width=.5in]{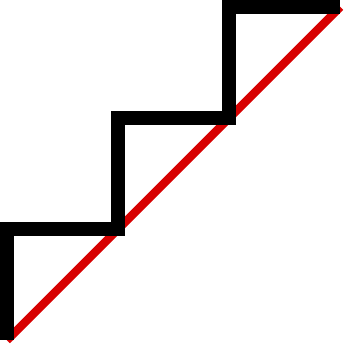}\hskip .2in
\includegraphics[width=.5in]{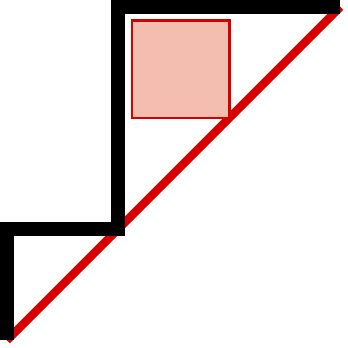}\hskip .2in
\includegraphics[width=.5in]{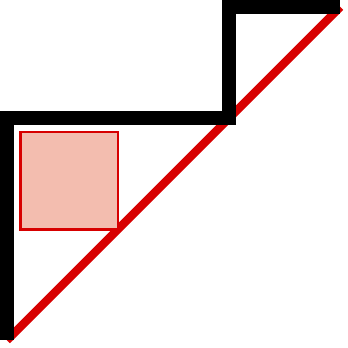}\hskip .2in
\includegraphics[width=.5in]{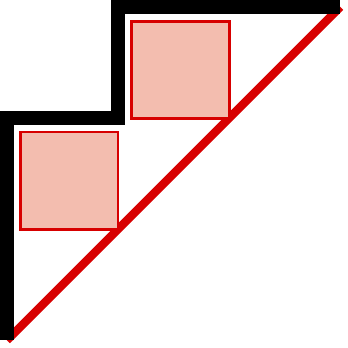}\hskip .2in
\includegraphics[width=.5in]{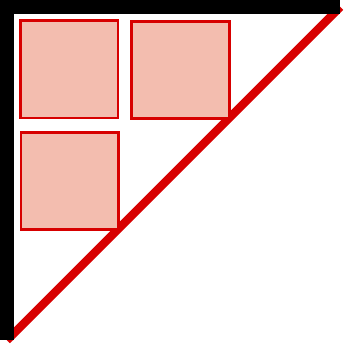}
\end{center}
\end{example}

We will make use of the Newton element 
\begin{align}
\label{newton}
\Omega[X] = \sum_{\lambda} p_\lambda[X]/z_\lambda 
=\sum_{m \geq 0} h_{m}[X], 
\end{align}
where we have the identities
\begin{align*}
\Omega[X + Y] &= \Omega[X]\Omega[Y]\\
\Omega[X - Y] &= \Omega[X]/\Omega[Y]\\
\Omega[x_1 + \epsilon x_2 - x_3 - \epsilon x_4] &= \frac{(1-x_3)(1+x_4)}{(1-x_1)(1+x_2)}\\
\Omega[XY] &= \sum_{\la} s_{\la}[X] s_\la[Y] = \sum_{\la} p_\la[X] p_\la[Y]/z_\la~.
\end{align*}

%*%Jing 
Jing \cite{Jin91} introduced  a family of 
operators $\Hop_m$ indexed by $m \in \ZZ$ using the following formal power
series in the parameter $z$,
\begin{align}\label{hopdef}
\Hop(z) P[X] &= \sum_{m \in \ZZ} z^m \Hop_m P[X] 
:= P\left[ X - \frac{1-q}{z} \right] \Omega[zX]\\
&=\sum_{m \in \mathbb Z} z^m 
 \sum_{r \geq 0} (-1)^r h_{m+r}[X] e_r[(1-q)X]^\perp P[X]~.\nonumber
\end{align}
He proved that these operators create the Hall-Littlewood polynomials
by adding rows.
\begin{prop} \cite{Jin91}  
%For $m \in {\mathbb Z}$, 
%let $\Hop_m P[X] = P\left[X - \frac{1-q}{z}\right] \Omega[zX] \coeff_{z^m}$
%where $\coeff_{z^m}$ represents the operation of taking the coefficient of $z^m$ in
%the expression to the left of the symbol.  This may also be expressed as
%\begin{equation}\label{hopdef}
%\Hop(z) P[X] = P\left[ X - \frac{1-q}{z} \right] \Omega[zX] = \sum_{m \in \ZZ} z^m \Hop_m P[X]~.
%\end{equation}
For any partition $\lambda=(\lambda_1,\ldots,\lambda_\ell)$,
$$Q'_\la[X;q] = \Hop_{\la_1} \Hop_{\la_2} \cdots \Hop_{\la_\ell}(1)~.$$
\end{prop}

\section{Two families of Hall-Littlewood symmetric functions}

Our primary focus is the study of two families of symmetric functions
and the combinatorics surrounding them.  These functions arise from
the following operators $\Bop_m$ and $\Cop_m$, 
closely related to Jing's $\Hop_m$ operators from equation \eqref{hopdef}:
\begin{align}\label{bopdef}
\Bop(z) P[X] &= \sum_{m \in {\mathbb Z}} z^m \Bop_m P[X] 
:= P\left[X + \epsilon \frac{(1-q)}{z}\right]\Omega[-\epsilon zX]\\
 &= \sum_{m \in {\mathbb Z}} z^m \sum_{r \geq 0} (-1)^r e_{m+r}[X] h_r[X(1-q)]^\perp P[X] \nonumber
\end{align}
and
\begin{align}\label{copdef}
\Cop(z) P[X]&=
 \sum_{m \in {\mathbb Z}} z^m \Cop_m P[X] 
:= -q~ P\left[X + \epsilon\frac{(1-q)}{z}\right]\Omega[\epsilon(z/q) X]\\
 &= \sum_{m \in {\mathbb Z}} (-1/q)^{m-1} z^m \sum_{r \geq 0} q^{-r} h_{m+r}[X] h_r[X(1-q)]^\perp P[X] \nonumber
\,.
\end{align}
The symmetric functions of particular interest here are those
defined, for any composition $\alpha$, by setting
$$
B_\alpha[X;q] = \Bop_{\alpha_{\ell(\alpha)}} 
\Bop_{\alpha_{\ell(\alpha)-1}} \cdots 
\Bop_{\alpha_1}(1)
$$
and
$$
C_\alpha[X;q] = \Cop_{\alpha_1} \Cop_{\alpha_2} 
\cdots \Cop_{\alpha_{\ell(\alpha)}}(1)~.
$$
Note that the operators generating $B_\alpha$ and $C_\alpha$ 
are both indexed by the parts of $\alpha$, but are
applied in reverse order with respect to one another.
This is done so that the associated
combinatorial and algebraic identities are more uniform.

These operators are related by way of the equation:
%\begin{align}
%-q \omega \Bop^{q\rightarrow1/q}(-z/q)\omega &= \Cop(z)
%\end{align}
%\begin{align}
%\omega \Bop(z)\omega &= \Hop(z)
%\end{align}
\begin{equation} \label{HandBrel}
\Bop(z) =  \omega \Hop( z ) \omega
\end{equation}
and
\begin{equation} \label{HandCrel}
\Cop(z) = (-q) \Hop^{q \rightarrow 1/q}(-z/q),
\end{equation}
%(note: replacing $z$ with $-z$ in the right hand side of \eqref{Bhopdef} 
%corresponds to replacing $z$ with $\epsilon z$ in the right hand side of the equation)
or equivalently 
\begin{equation}\label{operels}
\Cop_m = (-1/q)^{m-1} \Hop_m^{q \rightarrow 1/q} = (-1/q)^{m-1} \omega \Bop_m^{q \rightarrow 1/q} \omega~.
\end{equation}
Thus the functions themselves are related as:
\begin{equation} \label{BCrelation}
Q'_\lambda[X;q] = \omega B_{\rev{\lambda}}[X;q]
= (-q)^{ \ell(\lambda)-|\lambda|} 
C_\lambda[X;1/q] ~.
\end{equation}

The Jing operators create Hall-Littlewood polynomials indexed by partitions 
which form a basis for the symmetric function ring.  The $C_\alpha$ and $B_\alpha$
symmetric functions are 
indexed by compositions and are not linearly independent.  The 
equations above detail how Hall-Littlewood symmetric functions 
are included in these families and therefore, Schur positive expansions 
of the $C_\alpha$ and $B_\alpha$ hold in certain cases.  However,
they are not Schur positive in complete generality.  The smallest 
examples that are not uniformly Schur positive or Schur negative 
are $B_{(3,1)}[X;q]$ and $C_{(1,3)}[X;q]$.

\begin{example}  Below is a table of the symmetric functions $B_\alpha[X;q]$
and $C_\alpha[X;q]$ for $\alpha \models 4$.  Notice that both $B_{(3,1)}[X;q]$
and $C_{(1,3)}[X;q]$ have mixed signs in their coefficients.
\begin{equation*}
\begin{bmatrix}
B_{(1,1,1,1)}[X;q]\\
B_{(1,1,2)}[X;q]\\
B_{(1,2,1)}[X;q]\\
B_{(2,1,1)}[X;q]\\
B_{(1,3)}[X;q]\\
B_{(2,2)}[X;q]\\
B_{(3,1)}[X;q]\\
B_{(4)}[X;q]
\end{bmatrix}
=
\begin{bmatrix}
q^6&q^{3}+q^{4}+q^{5}&q^{2}+q^{4}&q+q^{2}+q^{3}&1\\
q^3&q+q^{2}&q&1&0\\
q^4&q^{2}+q^{3}&q^{2}&q&0\\
q^5&q^{3}+q^{4}&q^{3}&q^{2}&0\\
q&1&0&0&0\\
q^2&q&1&0&0\\
q^3&q^{2}&q-1&0&0\\
1&0&0&0&0
\end{bmatrix}
\begin{bmatrix}
s_{(1,1,1,1)}[X;q]\\
s_{(2,1,1)}[X;q]\\
s_{(2,2)}[X;q]\\
s_{(3,1)}[X;q]\\
s_{(4)}[X;q]
\end{bmatrix}
\end{equation*}

\begin{equation*}
\begin{bmatrix}
C_{(1,1,1,1)}[X;q]\\
C_{(1,1,2)}[X;q]\\
C_{(1,2,1)}[X;q]\\
C_{(2,1,1)}[X;q]\\
C_{(1,3)}[X;q]\\
C_{(2,2)}[X;q]\\
C_{(3,1)}[X;q]\\
C_{(4)}[X;q]
\end{bmatrix}
=
\begin{bmatrix}
1&q^{-3}+q^{-2}+q^{-1}&q^{-4}+q^{-2}&q^{-5}+q^{-4}+q^{-3}&q^{-6}\\
0&-q^{-3}&-q^{-4}&-q^{-5}-q^{-4}&-q^{-6}\\
0&-q^{-2}&-q^{-3}&-q^{-4}-q^{-3}&-q^{-5}\\
0&-q^{-1}&-q^{-2}&-q^{-3}-q^{-2}&-q^{-1}\\
0&0&q^{-3}-q^{-2}&q^{-4}&q^{-5}\\
0&0&q^{-2}&q^{-3}&q^{-4}\\
0&0&0&q^{-2}&q^{-3}\\
0&0&0&0&q^{-3}
\end{bmatrix}
\begin{bmatrix}
s_{(1,1,1,1)}[X;q]\\
s_{(2,1,1)}[X;q]\\
s_{(2,2)}[X;q]\\
s_{(3,1)}[X;q]\\
s_{(4)}[X;q]
\end{bmatrix}
\end{equation*}
\end{example}

To manipulate these symmetric functions, we derive commutation relations 
between the symmetric function operators. Our first result enables us
to expand an element $C_\alpha[X;q]$, for any composition $\alpha$,
in terms of the $C_\lambda[X;q]$ indexed by partitions $\lambda$.

\begin{prop}\label{Copconjrel}
For $m, n \in \ZZ$, we have
\begin{equation}
\label{Copconjreleq}
q \Cop_m \Cop_n - \Cop_{m+1} \Cop_{n-1} = \Cop_n \Cop_m - q \Cop_{n-1} \Cop_{m+1}~.
\end{equation}
In particular for $m \in \ZZ$,
\begin{equation}
\Cop_{m} \Cop_{m+1} = \frac{1}{q} \Cop_{m+1} \Cop_m~.
\end{equation}
\end{prop}

\begin{proof}
We begin with the expressions for the $\Cop$-operators from equation \eqref{copdef}.
For ease of notation we shall 
use $h_m$ in place of $h_m[X]$ and $h_r^{q\perp}$ in place of the expression $h_r[X(1-q)]^\perp$.
We compute that
\begin{align}
h_r^{q\perp} (h_m P[X]) &= h_m[X + (1-q)z] P[X + (1-q)z] \coeff_{z^r}\nonumber\\
&=\sum_{i\geq0} h_{m-i} h_i[1-q] P[X + (1-q)z] \coeff_{z^{r-i}} \label{hrcom}\\
&= \sum_{i \geq 0} h_{m-i} h_i[1-q] h_{r-i}^{q\perp} P[X]\nonumber.
\end{align}
We also know that 
\begin{equation}\label{hr1mqeval}
h_r[1-q] =
\begin{cases}
0 &\hbox{ if }r<0\,,\\
1 &\hbox{ if }r=0\,,\\
1-q &\hbox{ if }r>0\,.
\end{cases}
\end{equation}
These two identities imply that,
\begin{align}
q \Cop_m \Cop_n
= &(-1/q)^{m+n-2} \sum_{i\geq0}\sum_{r \geq 0} \sum_{s \geq 0} 
q^{-r-s-i+1} h_{m+r+i} h_{n+s-i} h_i[1-q] h_{r}^{q\perp} h_s^{q\perp}\nonumber\\
= &(-1/q)^{m+n-2} \sum_{i\geq0}\sum_{r \geq 0} \sum_{s \geq 0} 
q^{-r-s-i+1} h_{m+r+i} h_{n+s-i} h_{r}^{q\perp} h_s^{q\perp}\nonumber\\
&- (-1/q)^{m+n-2} \sum_{i\geq1}\sum_{r \geq 0} \sum_{s \geq 0} 
q^{-r-s-i+2} h_{m+r+i} h_{n+s-i} h_{r}^{q\perp} h_s^{q\perp}\nonumber\\
= &(-1/q)^{m+n-2} \sum_{i\geq0}\sum_{r \geq 0} \sum_{s \geq 0} 
q^{-r-s-i+1} h_{m+r+i} h_{n+s-i} h_{r}^{q\perp} h_s^{q\perp}\nonumber\\
&- (-1/q)^{m+n-2} \sum_{i\geq0}\sum_{r \geq 0} \sum_{s \geq 0} 
q^{-r-s-i+1} h_{m+r+i+1} h_{n+s-i-1} h_{r}^{q\perp} h_s^{q\perp}\nonumber\\
= &(-1/q)^{m+n-2} \sum_{r \geq 0} \sum_{s \geq 0} \sum_{i\geq0}
q^{-r-s-i+1} s_{(m+r+i,n+s-i)} h_{r}^{q\perp} h_s^{q\perp}.
\label{firstterm}
\end{align}
where to arrive at \eqref{firstterm} we have introduced the Schur function
$s_{(a,b)} = h_a h_b - h_{a+1} h_{b-1}$.
Similarly,
\begin{align*}
\Cop_{m+1} \Cop_{n-1} &= (-1/q)^{m+n-2} \sum_{r \geq 0} \sum_{s \geq 0} \sum_{i\geq0}
q^{-r-s-i} s_{(m+r+i+1,n+s-i-1)} h_{r}^{q\perp} h_s^{q\perp}\,.
\end{align*}

{}From these identities, we find that all terms in the difference 
$q \Cop_m \Cop_n - \Cop_{m+1} \Cop_{n-1}$
cancel except the $i=0$ term in \eqref{firstterm}:
\begin{align*}
q \Cop_m \Cop_n - \Cop_{m+1} \Cop_{n-1} = 
&(-1/q)^{m+n-2} \sum_{r \geq 0} \sum_{s \geq 0}
q^{-r-s+1} s_{(m+r,n+s)} h_{r}^{q\perp} h_s^{q\perp}\\
&+ (-1/q)^{m+n-2} \sum_{r \geq 0} \sum_{s \geq 0} \sum_{i\geq1}
q^{-r-s-i+1} s_{(m+r+i,n+s-i)} h_{r}^{q\perp} h_s^{q\perp}\\
&- (-1/q)^{m+n-2} \sum_{r \geq 0} \sum_{s \geq 0} \sum_{i\geq0}
q^{-r-s-i} s_{(m+r+i+1,n+s-i-1)} h_{r}^{q\perp} h_s^{q\perp}\\
= 
&(-1/q)^{m+n-2} \sum_{r \geq 0} \sum_{s \geq 0}
q^{-r-s+1} s_{(m+r,n+s)} h_{r}^{q\perp} h_s^{q\perp}.
\end{align*}
We can then compute $\Cop_n \Cop_m - q \Cop_{n-1} \Cop_{m+1}$
from this by replacing $m \rightarrow n-1$ and $n \rightarrow m+1$.
In particular,
\begin{align*}
q \Cop_{n-1} \Cop_{m+1} - \Cop_{n} \Cop_{m} = 
&(-1/q)^{m+n-2} \sum_{r \geq 0} \sum_{s \geq 0}
q^{-r-s+1} s_{(n-1+r,m+1+s)} h_{r}^{q\perp} h_s^{q\perp}.
\end{align*}
The identity $-s_{(b-1,a+1)} = s_{(a,b)}$ and the commutation 
of $h^\perp$ then implies our claim.
%Note: Doesn't this just follow from the identity by putting n=m+1?  I don't
%think we need to include this paragraph:
%The same relation still holds when $n = m+1$, but in this case there 
%is simplification in
%the identity and $q \Cop_m \Cop_{m+1} - \Cop_{m+1} \Cop_{m} = \Cop_{m+1} \Cop_m - q \Cop_{m} \Cop_{m+1}$
%reduces to $q \Cop_m \Cop_{m+1} = \Cop_{m+1} \Cop_{m}$.
\end{proof}

An important consequence of this result is that if $\alpha$ is a
composition of length $\ell$, then $C_{\alpha}[X;q]$ can be written
as a linear combination of the $C_{\lambda}[X;q]$ where $\lambda$ are
partitions that also have length $\ell$.
\begin{example}
The symmetric function $C_{(1,3)}[X;q]$ can be expressed in terms of
$C_{(3,1)}[X;q]$ and $C_{(2,2)}[X;q]$ using this commutation relation
since $q \Cop_1 \Cop_3 = \Cop_2 \Cop_2 + \Cop_3 \Cop_1 - q \Cop_2 \Cop_2$.  
Consequently, $C_{(1,3)}[X;q] = (1/q -1)C_{(2,2)}[X;q] + 1/qC_{(3,1)}[X;q]$.
\end{example}

The relation of $\Cop$ to $\Bop$ given in \eqref{operels} 
enables us to derive an identity on $\Bop$ from
Theorem~\ref{Copconjrel}.
In particular, we simply apply $\omega$ to \eqref{Copconjreleq}
and replace $q$ by $1/q$.

\begin{cor}  For $m \in \ZZ$,
$$ \Bop_m \Bop_n - q\Bop_{m+1} \Bop_{n-1} 
= q \Bop_n \Bop_m - \Bop_{n-1} \Bop_{m+1}\,.$$
In particular, letting $n = m+1$ gives
$$ \Bop_m \Bop_{m+1} 
= q \Bop_{m+1} \Bop_m ~.$$
\end{cor}

In fact, we can also pin down commutation relations
between the $\Bop$ and $\Cop$ operators if
$m+n>0$ (note, the relation does not hold when
$m+n\leq 0$).

\begin{prop} \label{BCcommutation}
If $m + n > 0$, then
\begin{equation}
\Bop_n \Cop_m = q \Cop_m \Bop_n~.
\end{equation}
\end{prop}

\begin{proof}
We use identities \eqref{hrcom} and \eqref{hr1mqeval} and
compute an expression for $\Bop_m \Cop_n$:
\begin{align*}
\Bop_m \Cop_n = &(-1/q)^{n-1}\sum_{r \geq 0}\sum_{s \geq 0} \sum_{i \geq 0}
(-1)^r q^{-s} e_{m+r}  h_{n+s-i} h_i[1-q] h_{r-i}^{q\perp} h_s^{q\perp}\\
=&(-1/q)^{n-1} \sum_{r \geq 0}\sum_{s \geq 0} 
(-1)^{r} q^{-s} e_{m+r}  h_{n+s} h_{r}^{q\perp} h_s^{q\perp}\\
&+(-1/q)^{n-1} \sum_{i \geq 1} \sum_{r \geq 0}\sum_{s \geq 0} 
(-1)^{r+i} q^{-s} e_{m+r+i}  h_{n+s-i} (1-q) h_{r}^{q\perp} h_s^{q\perp}.
\end{align*}
Analogously, we also have the equations 
$h_r^{q\perp} e_m = \sum_{i \geq 0} e_{m-i} e_i[1-q] h_{r-i}^{q\perp}$
and $e_r[1-q] =
(-q)^{r-1}h_r[1-q]\hbox{ if }r>0$.  From this we derive a similar 
expression for $q \Cop_n \Bop_m$:
\begin{align*}
q \Cop_n \Bop_m 
%&= 
%q (-1/q)^{n-1} \sum_{s \geq 0}\sum_{r \geq 0}(-1)^r q^{-s} h_{n+s} h_s^{q\perp}
% e_{m+r}  h_r^{q\perp}\\
= 
&(-1/q)^{n-1} \sum_{s \geq 0}\sum_{r \geq 0}\sum_{i \geq 0}
(-1)^r q^{-s+1} h_{n+s}  
 e_{m+r-i} e_i[1-q] h_{s-i}^{q\perp} h_r^{q\perp}\\
= 
&(-1/q)^{n-1} \sum_{s \geq 0}\sum_{r \geq 0}
(-1)^r q^{-s+1} e_{m+r} h_{n+s}  
 h_{s}^{q\perp} h_r^{q\perp}\\
&+(-1/q)^{n-1} \sum_{i \geq 1} \sum_{s \geq 0}\sum_{r \geq 0}
(-1)^{r+i+1} q^{-s} e_{m+r-i} h_{n+s+i}  
 (1-q) h_{s}^{q\perp} h_r^{q\perp}.
\end{align*}

Their difference is 
\begin{align*}
\Bop_m\Cop_n  - q \Cop_n \Bop_m =
&(1-q)(-1/q)^{n-1} \sum_{s \geq 0}\sum_{r \geq 0}
q^{-s} (-1)^r e_{m+r} h_{n+s}  
 h_{r}^{q\perp} h_s^{q\perp}\\
&+(1-q)(-1/q)^{n-1} \sum_{r \geq 0} \sum_{s \geq 0}q^{-s}
\left(\sum_{i \geq 1} 
(-1)^{r+i} e_{m+r-i}  h_{n+s+i}\right) h_{r}^{q\perp} h_s^{q\perp}\\
&+(1-q)(-1/q)^{n-1} \sum_{r \geq 0} \sum_{s \geq 0}q^{-s}
\left(\sum_{i \geq 1} 
(-1)^{r+i} e_{m+r+i}  h_{n+s-i}\right) h_{r}^{q\perp} h_s^{q\perp}\,.
\end{align*}
In fact, the right hand side reduces to zero since
$m+n>0$ implies that for each $r,s\geq 0$, $m+n+r+s>0$ and
\begin{align*}
(-1)^r e_{m+r} h_{n+s} 
&+\sum_{i \geq 1} (-1)^{r+i} e_{m+r-i} h_{n+s+i}
+ \sum_{i\geq1} (-1)^{r+i} e_{m+r+i}  h_{n+s-i} \\
&= \sum_{i=-n-s}^{m+r} (-1)^{r+i} e_{m+r-i} h_{n+s+i}=0
\end{align*}
by the identity $\sum_{i=0}^d (-1)^i e_{d-i} h_i = 0$ for all $d>0$.
\end{proof}

Jing's operators generalize operators of Bernstein (see \cite{Macdonald}),
defined by
\begin{align}\label{Bhopdef}
\Sop(z) P[X] &= \sum_{m \in \ZZ} z^m \Sop_m P[X] = P\left[ X - \frac{1}{z} \right] \Omega[zX]\\
&= \sum_{m \in \ZZ} z^m \sum_{r \geq 0} (-1)^r h_{m+r}[X] e_r[X]^\perp P[X]~.\nonumber
\end{align}
These are creation operators for the Schur functions since
$\Sop_{\la_1} \Sop_{\la_2} \cdots \Sop_{\la_\ell}(1) = s_{\lambda}[X]$,
and they satisfy the commutation relation $\Sop_m \Sop_n = - \Sop_{n-1} \Sop_{m+1}$.
We can write the Schur creation operators in terms of the $\Cop_a$ operators, which will
help us in \S~\ref{Sym} to write Schur functions in terms of the $C_{\alpha}$.
\begin{prop}
For $m \in \ZZ$,
\begin{equation}\label{SopexpCop}
\Sop_m = (-q)^{m-1} \sum_{i \geq 0} \Cop_{m+i} e_i^\perp .
\end{equation}
\end{prop}
\begin{proof}
We will use the identity $h_r[(1-q)X] = \sum_{j\geq0} h_j[X] h_{r-j}[-qX] = \sum_{j\geq0} (-q)^{r-j} h_j[X] e_{r-j}[X]$ and calculate directly,
\begin{align*}
(-q)^{m-1} \sum_{i \geq 0} \Cop_{m+i} e_i^\perp &=
(-q)^{m-1} \sum_{i \geq 0} (-1/q)^{m+i-1} \sum_{r \geq 0} q^{-r} h_{m+i +r}[X] h_r[X(1-q)]^\perp e_i^\perp\\
&=
\sum_{r \geq 0} \sum_{i \geq 0} \sum_{j\geq0} (-1/q)^{i-r+j} q^{-r} h_{m+i +r}[X] h_j^\perp 
e_{r-j}^\perp e_i^\perp\\
&=
\sum_{d \geq 0} \sum_{j=0}^d \sum_{r \geq 0}  (-1)^{d-r} q^{-d} h_{m+d+r-j}[X] h_j^\perp 
e_{r-j}^\perp e_{d-j}^\perp\\
&=
\sum_{d \geq 0} \sum_{j=0}^d \sum_{r \geq 0} (-1)^{d-r-j} q^{-d} h_{m+d+r}[X] h_j^\perp e_{r}^\perp 
e_{d-j}^\perp\\
&=
\sum_{r \geq 0} (-1)^{-r} h_{m+r}[X] e_{r}^\perp = \Sop_m~
\end{align*}
where the last equality follows because $\sum_{j=0}^d (-1)^j h_j e_{d-j} = 0$ for
all $d > 0$ so the remaining sum is only the part where $d=0$.
\end{proof}
%
%We first use \eqref{ekperp} to compute that
%\begin{align*}
%\Sop_m P[X] &= P \left[ X - \frac{1}{z} \right] \Omega[zX] \coeff_{z^m}\\
%&= (-q)^{m-1} (-q)^{-m+1} P\left[X - \frac{(1-1/q)}{z} - \frac{1}{qz}\right]\Omega[zX]\coeff_{z^m}\\
%&= (-q)^{m-1}  \sum_{i \geq 0} (-q)^{-m-i+1}\left(\frac{1}{z}\right)^i e_i^\perp P\left[X - \frac{(1-1/q)}{z} \right]\Omega[zX]\coeff_{z^m}\,.
%\end{align*}
%The result then follows by noting that
%$\Cop_a e_i^\perp = e_i^\perp \Cop_a$ and
%alternatively writing the action of $\Cop$ as

\begin{remark}
In the reference \cite{GXZ10}, the operator $\Cop_a$ is presented
in a slightly different but equivalent expression.
We note that a series $f(z) = \sum_{n \in \ZZ} f_n z^n$ has
the property that $f(\epsilon z/q) \coeff_{z^a} = (-1/q)^a f(z) \coeff_{z^a}$.
For this reason,
\begin{align*}
\mathbb C _a P[X] &= -q P \left [ X + \epsilon \frac {(1-q)}{z} \right ] 
\Omega [ \epsilon (z/q) X ] \coeff_{z^a} \\
&= -q (-1/q)^a P \left [ X + \frac {(1-q)}{q z} \right ] 
\Omega [ z X ] \coeff_{z^a} \\
&= (-1/q)^{a-1} P \left [ X - \frac {1-1/q}{z} \right ] \Omega [zX] \coeff_{z^a}.
\end{align*}
\end{remark}

\section{The combinatorics of $\nabla$ applied to Hall-Littlewood polynomials} \label{Nab}

%The primary motivation behind the exploration of
%$B_\alpha[X;q]$ and $C_\alpha[X;q]$ is our discovery that 
%the combinatorics of Dyck paths can naturally be refined by 
%the study of $\nabla(B_\alpha[X;q])$ and $\nabla(C_\alpha[X;q])$. As a
%consequence we have used this as a guide for the discovery
%of symmetric function identities.

Recall that in the special case that $\alpha$ is a partition, 
$B_\alpha[X;q]$ and $C_\alpha[X;q]$ are closely related to the 
Hall-Littlewood symmetric functions.  It was conjectured in 
\cite[Conjecture II and III]{BGHT99} (partially attributed 
to A. Lascoux) that applying $\nabla$ to a Hall-Littlewood polynomial
produces a Schur positive function.  Our main discovery is that
including all compositions $\alpha$ in the study of 
$\nabla(B_\alpha[X;q])$ and $\nabla(C_\alpha[X;q])$
leads to a natural refinement for the combinatorics of Dyck paths.
Moreover, our combinatorial exploration led
us to discover new symmetric function identities.

One useful tool in the exploration of the operator $\nabla$ is the
fact that $\nabla^{q=1}$ is a multiplicative operator.
Since we can deduce from the operator definitions of
our symmetric functions that $B_\alpha[X;1] = e_\alpha[X]$ 
and $C_\alpha[X;1] = h_\alpha[X]$,
we have
$$\nabla^{q=1}( B_\alpha[X;1] ) = \nabla^{q=1}( e_{\alpha_\ell}[X] )\nabla^{q=1}( e_{\alpha_{\ell-1}[X]} )\cdots \nabla^{q=1}( e_{\alpha_1}[X] )$$
and
$$
\nabla^{q=1}( C_\alpha[X;1] ) = \nabla^{q=1}( h_{\alpha_1}[X] )\nabla^{q=1}
( h_{\alpha_2}[X] )\cdots \nabla^{q=1}( h_{\alpha_\ell}[X] )~.
$$
From this we can deduce the coefficient of $e_n[X]$.
In 
particular, the coefficient of $e_n[X]$ in $\nabla( e_n[X] )$ and 
in $\nabla( h_n[X] )$ is the $q,t$-Catalan numbers $C_n(q,t)$ and 
$C_{n-1}(q,t)$, respectively.
Thus, the coefficient of $e_n[X]$ in $\nabla^{q=1}( B_\alpha[X;1] )$
and in $\nabla^{q=1}( C_\alpha[X;1] )$ is $\prod_i C_{\alpha_i}(1,t)$ and
$\prod_i C_{\alpha_i-1}(1,t)$ respectively.  
The combinatorial interpretation for $C_n(1,t)$ then
gives combinatorial meaning to these coefficients.
Namely, $\left< e_n[X], \nabla^{q=1}( B_\alpha[X;1] )\right>$ is
the $t$-enumeration of Dyck paths (with weight $t$ raised to the \area) 
which lie below the staircase consisting of $\alpha_1$ steps
up and over, $\alpha_2$ steps up and over, etc.  
and $\left< e_n[X], \nabla^{q=1}( C_\alpha[X;1] )\right>$ is a $t$-enumeration of
Dyck paths which touch the diagonal only in rows $1$, $1+\alpha_1$, $1+\alpha _1+\alpha_2$ steps, etc.

\begin{prop}  For $\alpha$ a composition of $n$,
\begin{equation}\label{qeq1nablaBalpha}
\left< \nabla^{q=1}( B_\alpha[X;1] ), e_n[X] \right> = \sum_{D \leq DP(\alpha)} t^{\area(D)}
\end{equation}
and
\begin{equation}\label{qeq1nablaCalpha}
\left< \nabla^{q=1}( C_\alpha[X;1] ), e_n[X] \right> = \sum_{\touch(D) = \alpha} t^{\area(D)}.
\end{equation}
\end{prop}

Remarkably, we have empirical evidence to suggest that
in general, there is a combinatorial interpretation for the 
coefficient of $e_n[X]$ in $\nabla(B_\alpha[X;q])$ and 
in $\nabla(C_\alpha[X;q])$ that naturally
generalizes the beautiful combinatorics of the $q,t$-Catalan.

\begin{conj} \label{belowCI}
For $\alpha \models n$,
\begin{equation}\label{nBa}
\left< \nabla( B_\alpha[X;q] ), e_n[X] \right>  
= \sum_{D \leq DP(\alpha)}  t^{\area(D)} q^{\dinv(D)+\doff_\alpha(D)}~.
\end{equation}
\end{conj}

\begin{conj} \label{touchCI}
For $\alpha \models n$,
\begin{equation}\label{nCa}
\left< \nabla(C_\alpha[X;q]), e_n[X] \right> 
= \sum_{\touch(D) = \alpha} t^{\area(D)} q^{\dinv(D)}. 
\end{equation}
\end{conj}

Our work was inspired by the work of %*%Bergeron, Descouens, and Zabrocki
\cite{BDZ10} where they considered coefficients
$\nabla( B_{\rev{\lambda}}[X;q])$ for $\lambda$ a hook partition (since
for that case $\doff_{\rev{\lambda}}(D) = 0$).
The innovation in these identities is to consider symmetric functions 
indexed by compositions which allowed us to conjecture the action of
$\nabla$ on a spanning set of the symmetric functions.

More generally, we have conjectures for the expansion of
$\nabla( B_{\alpha}[X;q] )$ and
$\nabla( C_{\alpha}[X;q] )$ into monomials.  

\begin{conj}\label{combinterpBa}
\begin{equation}\label{nBa2}
\nabla( B_\alpha[X;q] )  = \sum_{D \leq DP(\alpha)} 
\sum_{w \in \WP_{D}} t^{\area(w)} 
q^{\dinv(w)+\doff_\alpha(D)} x^w.
\end{equation}
\end{conj}

\begin{conj}\label{combinterpCa}
\begin{equation}\label{nCa2}
\nabla(C_\alpha[X;q]) = \sum_{\touch(D) = \alpha} \sum_{w \in \WP_{D}} t^{\area(D)} q^{\dinv(w)} x^w.
\end{equation}
\end{conj}

By the arguments in \cite{HHLRU05}
(see also \cite[p.99]{Hag08}) Conjecture \ref{combinterpBa} and
\ref{combinterpCa} imply
Conjecture \ref{belowCI} and \ref{touchCI}.
The case $\alpha = (n)$ of \eqref{nBa2} reduces to the shuffle conjecture since $B_{(n)}[X;q] = e_n[X]$.
Also, because of the expansion of $s_{(n-k,1^k)}$ in Proposition
\ref{sn1nmkexpansion} in the next section, \eqref{nCa2} implies the special 
case of the Loehr-Warrington conjecture \cite[Conjecture 3]{LoWa08} involving the 
action of $\nabla$ on the Schur function $s_{(n-k,1^k)}$.  

We will prove in the next section 
that Conjecture \ref{combinterpBa} and \ref{combinterpCa} are equivalent to
each other (and by consequence Conjecture \ref{belowCI} and \ref{touchCI} are
equivalent as well).  In work building on 
our results here, %*%Garsia, Xin and Zabrocki 
\cite{GXZ10} with contributions 
from %*%Hicks 
\cite{Hic10} proved Conjecture \ref{touchCI}.

\section{Symmetric Function Identities}
\label{Sym}

%The comparison of the combinatorics of our generalized
%$q,t$-Catalans to their algebraic counterparts has brought
%to light a number of symmetric function identities involving
%$B_\alpha[X;q]$ and $C_\alpha[X;q]$.  Recall that in the
The exploration of $q,t$-Catalans led %*%Garsia and Haglund 
\cite{GaHa01}
to the special symmetric function elements 
$E_{n,k}[X;q]$, defined by the algebraic identity
$$
e_n\left[X\frac{1-z}{1-q} \right] = \sum_{k=1}^n \frac{(z;q)_k}{(q;q)_k} E_{n,k}[X;q]
$$
where $(z;q)_k = (1-z)(1-qz)\cdots(1-q^{k-1}z)$.
These elements play a fundamental role in the proof
that the $q,t$-Catalan polynomial is the $q,t$-enumeration 
of Dyck paths as given in \eqref{CIqtcatalan}.
Namely, the proof follows by showing that
\begin{equation}
\label{recur}
\langle\nabla E_{n,k}[X;q], e_n[X]\rangle=
q^{k\choose 2} t^{n-k} \sum_{r=0}^{n-k}
\begin{bmatrix}
r+k-1 \\
r
\end{bmatrix}_q
\left<\nabla(E_{n-k,r}[X;q]),e_{n-r}[X] \right>\,,
\end{equation}
where $\begin{bmatrix}n\\k\end{bmatrix}_q = \frac{(q;q)_n}{(q;q)_k (q;q)_{n-k}}$,
and the combinatorial interpretation for
$\langle\nabla E_{n,k}[X;q], e_n[X]\rangle$ in terms of Dyck paths
satisfies the same recurrence.

In particular, the coefficient of $e_n[X]$ in $\nabla( E_{n,k}[X;q] )$ 
$q,t$-enumerates the Dyck paths which touch the diagonal $k$ times.  
{}From this, Conjecture~\ref{touchCI} leads us to expect that
$$
\left< e_n[X], \nabla( E_{n,k}[X;q] ) \right> = 
\sum_{\alpha \models n, \ell(\alpha) = k} \left< e_n[X], 
\nabla( C_{\alpha}[X;q] ) \right>\,.
$$
In fact, we have discovered much more generally that
\begin{equation}
E_{n,k}[X;q]=\sum_{\alpha \models n\atop\ell(\alpha)=k}
C_\alpha[X;q]\,.
\end{equation}
This section is devoted to proving this surprising result, which 
suggests that the $C_\alpha[X;q]$ are the building blocks in 
$q,t$-Catalan theory.

\begin{remark}
A key step in the proof of our Conjectures \ref{belowCI} 
and \ref{touchCI} relies on extending the recurrence \eqref{recur} 
to involve Dyck paths which touch the diagonal at certain points.
This is carried out in \cite{GXZ10}.
\end{remark}

Our point of departure is to give a simple expression for $e_n[X]$ in
terms of the Hall-Littlewood symmetric functions $C_\alpha[X;q]$.

\begin{prop} \label{enexpansioninCa}
\begin{equation}\label{enexpansion}
e_n[X] = \sum_{\alpha \models n} C_\alpha[X;q].
\end{equation}
\end{prop}
\begin{proof}
Assume by induction on $n$ that equation \eqref{enexpansion} holds (the base cases of
$n=0$ and $1$ are easily verified).  Since the operator $\Sop_m$
is a creation operator for the Schur functions, by \eqref{SopexpCop} we have
\begin{align*}
e_n[X] &= s_{(1^n)}[X]
= \Sop_1( s_{(1^{n-1})}[X])
= \sum_{i=0}^{n-1} \Cop_{1+i} s_{(1^{n-i-1})}[X]
\end{align*}
Which, by induction, gives
\begin{align*}
e_n[X]=
\sum_{i=0}^{n-1} \sum_{\alpha\models n-i-1} \Cop_{i+1} C_{\alpha}[X;q]
= \sum_{\alpha\models n} C_\alpha[X;q]~.
\end{align*}
\end{proof}

Proposition~\ref{enexpansioninCa} can be stated in a 
more general form, suggesting that any Schur function 
may expand nicely  in terms of our
Hall-Littlewood spanning set.

\begin{prop}\label{sn1nmkexpansion} For $0 \leq k < n$,
$$s_{(n-k,1^k)}[X] = (-q)^{n-k-1}\sum_{\substack{\alpha \models n\\\alpha_1 \geq n-k}} 
C_\alpha[X;q]~.$$
\end{prop}

\begin{proof}
Again using that $\Sop$ is a Schur function creation operator,
by \eqref{SopexpCop} we have
\begin{align*}
s_{(n-k,1^k)}[X] 
&= \Sop_{n-k}( s_{1^k}[X] )
= (-q)^{n-k-1} \sum_{i=0}^{k} \Cop_{n-k+i} \left( s_{(1^{k-i})}[X] \right)
\,.
\end{align*}
The previous proposition then implies
\begin{align*}
s_{(n-k,1^k)}[X] 
&= (-q)^{n-k-1} \sum_{i=0}^{k} \Cop_{n-k+i} \left(\sum_{\alpha \models k-i} 
C_{\alpha}[X;q]\right)\\
&= (-q)^{n-k-1} \sum_{i=0}^{k} \sum_{\alpha \models k-i} C_{(n-k+i,\alpha)}[X;q]
= (-q)^{n-k-1} \sum_{\substack{\alpha \models n\\
\alpha_1\geq n-k}} C_{\alpha}[X;q]~.
\end{align*}
\end{proof}

We are now in the position to prove that $E_{n,k}[X;q]$ can be
decomposed canonically in terms of the $C_\alpha[X;q]$.

\begin{prop}\label{Enkexpansion}
For $0\leq k<n$,
\begin{equation}\label{Enkformula}
%E_{n,k}[X;q] = \sum_{\begin{subarray}{c}\mu \vdash n\\\ell(\mu) = k\end{subarray}} 
E_{n,k}[X;q] = \sum_{\substack{\mu \vdash n\\\ell(\mu) = k}} 
q^{-n(\mu)-k + M(\mu)}
\bchoose{k}{m(\mu)}_q
C_\mu[X;q], 
%q^{\sum _{i=1}^n \left( \substack{m_i(\mu)+1\\2}\right)}
\end{equation}
where $M(\mu)=\sum_{i=1}^n {m_i(\mu)+1 \choose 2}$
and $\bchoose{k}{m(\mu)}_q=
(q;q)_k/\prod_{i=1}^n (q;q)_{m_i(\mu)}$.
\end{prop}

\begin{proof}
Recall the expansion of the elementary symmetric functions in the Macdonald basis
is given by (see \cite{GaHa01})
$$
e_n\left[X\frac{1-z}{1-q}\right] = \sum_{\mu \vdash n} 
\frac{\Ht_\mu[X;q,t] \Ht_\mu[(1-z)(1-t);q,t]}{{\tilde h}_\mu(q,t) {\tilde h}'_\mu(q,t)},
$$
where ${\tilde h}_\mu(q,t) = \prod_{c \in \mu} (q^{a(c)} - t^{l(c) +1})$ and
${\tilde h}'_\mu(q,t) = \prod_{c \in \mu}(t^{l(c)} - q^{a(c) +1})$.  
When $t=0$, these expressions reduce to
$${\tilde h}_\mu(q,0) = \prod_{c \in \mu} (q^{a(c)} - 0^{l(c) +1}) = q^{n(\mu')}$$
and
\begin{align*}
{\tilde h}'_\mu(q,0) &= \prod_{\substack{c \in \mu\\l(c)=0}}(1 - q^{a(c) +1})
\prod_{\substack{c \in \mu\\l(c)\neq 0}}(- q^{a(c) +1}) \\
&=(-1)^{n-\mu_1} q^{n+n(\mu')-M(\mu ^{\prime})} \prod_{i=1}^n (q;q)_{m_i(\mu')}.
\end{align*}
Therefore since $\Ht_\mu[X;q,t] = \Ht_{\mu'}[X;t,q]$, when we set $t=0$ everywhere we have
\begin{equation}
\label{ent0}
e_n\left[X\frac{1-z}{1-q}\right] = \sum_{\mu \vdash n} 
(-1)^{n-\mu_1} q^{-n-2n(\mu')+M(\mu ^{\prime})} 
\frac{\Ht_{\mu'}[X;0,q] \Ht_{\mu}[(1-z);q,0]}
{\prod_{i=1}^n (q;q)_{m_i(\mu')}}\,.
\end{equation}
Now the evaluation
\begin{equation}
\Ht_\mu[(1-z);q,t] = \prod_{c \in \mu} (1-z t^{l'(c)} q^{a'(c)})
\end{equation}
also yields
$$\Ht_\mu[(1-z);q,0] = (z;q)_{\mu_1}~.$$
Thus, replacing $\mu$ by $\mu ^{\prime}$ in \eqref{ent0},
and thereby exchanging $\mu_1$ and $\ell(\mu)$, gives
\begin{equation}
e_n\left[X\frac{1-z}{1-q}\right] = \sum_{\mu \vdash n} 
(-1)^{n-\ell(\mu)} q^{-n-2n(\mu)+M(\mu)}
\frac{\Ht_{\mu}[X;0,q] (z;q)_{\ell(\mu)}}
{\prod_{i=1}^n (q;q)_{m_i(\mu)}}\,.
\end{equation}
Since
$\Ht_\mu[X;0,q] = (-1)^{n-\ell(\mu)} q^{n(\mu) + n - \ell(\mu)} C_\mu[X;q]$, we
have
$$
e_n\left[X\frac{1-z}{1-q}\right] = \sum_{k=1}^n \frac{(z;q)_k}{(q;q)_k} 
\sum_{\substack{\mu \vdash n\\\ell(\mu)=k}} 
q^{-k-n(\mu)+M(\mu)}
\frac{C_\mu[X;q] (q;q)_{k}}
{\prod_{i=1}^n (q;q)_{m_i(\mu)}}\,,
$$
which implies our claim.
\end{proof}
The $q$-binomial coefficients that appear in equation \eqref{Enkformula}
suggest that there is a relation between the terms of the
$C_\la[X;q]$ basis and subsets of a $k$ element set.  It turns out 
that Proposition \ref{Enkexpansion} can be more cleanly written 
over compositions using a different expansion.

\begin{cor} \label{cleanEnkformula}
For $0\leq k<n$,
\begin{equation}
E_{n,k}[X;q] = \sum_{\substack{\alpha\models n\\ \ell(\alpha)=k}} C_\alpha[X;q]\,.
\end{equation}
\end{cor}
\begin{proof}
Using the straightening relations of the $\Cop_m$ operators,
if $\alpha$ is a composition of $n$ and $\la$ is a partition
of $n$ such that  $\ell(\la) \neq \ell(\alpha)$, then
$$C_\alpha[X;q] \coeff_{C_\la[X;q]} = 0~.$$
Now  for $\ell(\la) = k$, by Proposition \ref{Enkexpansion} and the fact that $e_n[X] = \sum_{k=1}^n E_{n,k}[X;q]$, 
\begin{align*}
E_{n,k}[X;q] \coeff_{C_\la[X;q]} &= e_n[X] \coeff_{C_\la[X;q]}\\
&= \sum_{\alpha \models n} C_\alpha[X;q] \coeff_{C_\la[X;q]}\\
&= \sum_{\substack{ \alpha \models n\\ \ell(\alpha)= k}} C_\alpha[X;q] \coeff_{C_\la[X;q]}.
\end{align*}
Furthermore if $\ell(\la) \neq k$ then
$$
E_{n,k}[X;q] \coeff_{C_\la[X;q]} = 0 = 
\sum_{\substack{ \alpha \models n\\ \ell(\alpha)= k}} C_\alpha[X;q] \coeff_{C_\la[X;q]}.  
$$
Since the functions $C_\la[X;q]$ are a basis this implies that
$E_{n,k}[X;q] = \sum_{\substack{ \alpha \models n\\ \ell(\alpha)= k}} C_\alpha[X;q]$.
\end{proof}

%These expressions imply why Conjectures \ref{touchCI} and \ref{combinterpCa} 
%are likely to hold given the results of \cite{GaHa01, GaHa02, Hag08, HHLRU05} 
%which provide a conjectured interpretation for 
%$\left< e_n[X], \nabla( E_{n,k}[X;q] ) \right>$
%and $\nabla( E_{n,k}[X;q] )$.  This shows that these new conjectures 
%refine previously
%known results and conjectures because there is not in general an expansion of
%$C_\alpha[X;q]$ in terms of other elements where an application of $\nabla$ is
%understood.

We have seen in \eqref{operels} that $\Cop$ is naturally
related to $\Bop$.  Here we pin down the relationship between 
the symmetric functions $B_\alpha[X;q]$ and $C_\alpha[X;q]$.
A by-product of this identity is that Conjecture~\ref{combinterpCa} 
implies Conjecture \ref{combinterpBa}.  

\def\mleq{{\leq\!\cdot}}

\begin{thrm} 
\label{BatoCaexpansion}
For $n \geq 0$ and any composition $\alpha \models n$,
\begin{equation} \label{BatoCa}
B_{\alpha}[X;q] = \sum_{\beta \leq \alpha} q^{\doff_\alpha(DP(\beta))} 
C_\beta[X;q]\,.
\end{equation}
\end{thrm}

\begin{proof}
We show this result by induction on the number of parts of $\alpha$.  The
base case follows since $\Bop_m(1) = e_m[X]$ which is equal to 
$\sum_{\gamma\models m} C_\gamma[X;q]$ by Proposition~\ref{enexpansioninCa}.
Assume by induction that \eqref{BatoCa} holds for a composition 
$\alpha$ of length $\ell$ and consider a composition $(m, \alpha)$.
We then have 
\begin{align}
\label{Balpham}
B_{(\alpha,m)}[X;q] &= \Bop_m( B_\alpha[X;q])
=\sum_{\beta \leq \alpha} q^{\doff_\alpha(DP(\beta))} \Bop_m( C_\beta[X;q] )\,.
%  \\
%&=\sum_{\beta \leq \alpha} q^{\doff_\alpha(DP(\beta))+\ell(\beta)} 
%\sum_{\gamma \models m} C_{(\beta,\gamma)}[X;q].
\end{align}

Now consider
$\Bop_m( C_\beta[X;q] ) = \Bop_m \circ \Cop_{\beta_1}\circ \Cop_{\beta_2} 
\circ \cdots \circ\Cop_{\beta_{\ell(\beta)}}(1)$.  The commutation 
relation between the $\Cop_n$ and $\Bop_m$ from Theorem \ref{BCcommutation} 
implies
\begin{align*}
\Bop_m( C_\beta[X;q] ) 
&= q^{\ell(\beta)} \Cop_{\beta_1}\circ \Cop_{\beta_2} 
\circ \cdots \circ\Cop_{\beta_{\ell(\beta)}} \circ \Bop_m(1)\,.
%&= q^{\ell(\alpha)} \Cop_{\alpha_1}\circ \Cop_{\alpha_2} 
%\circ \cdots \circ\Cop_{\alpha_{\ell(\alpha)}} (e_m[X])\\
\end{align*}
By Proposition~\ref{enexpansioninCa}, we then have
\begin{align}
\Bop_m( C_\beta[X;q] ) 
&= q^{\ell(\beta)} \Cop_{\beta_1}\circ \Cop_{\beta_2} 
\circ \cdots \circ\Cop_{\beta_{\ell(\beta)}} 
( \sum_{\gamma \models m} C_\gamma[X;q] )\\
&= q^{\ell(\beta)}  \sum_{\gamma \models m} C_{(\beta,\gamma)}[X;q] )\,.
\label{BmC}
\end{align}

Putting \eqref{BmC} into \eqref{Balpham}, we thus find that
\begin{align*}
B_{(\alpha,m)}[X;q] 
&=\sum_{\beta \leq \alpha} q^{\doff_\alpha(DP(\beta))+\ell(\beta)} 
\sum_{\gamma \models m} C_{(\beta,\gamma)}[X;q]\,.
\end{align*}
For each term in the sum, the composition $(\beta,\gamma)$ is finer 
than the composition $(\alpha, m)$.  Moreover, if we let $r_i$ be the 
number of times that $DP(\beta)$ touches the diagonal below the
$i^{th}$ bump of the Dyck path $DP(\alpha)$, then
$\doff_{(\alpha,m)}( DP(\beta,\gamma) ) = 
\sum_{i=1}^{\ell(\alpha)} r_i (\ell(\alpha)+1-i)
= \sum_{i=1}^{\ell(\alpha)} r_i (\ell(\alpha)-i) + \sum_{i=1}^{\ell(\alpha)} r_i
= \doff_{\alpha}(DP) + \ell(\beta)$.
\end{proof}

To prove that that Conjectures~\ref{combinterpCa} and \ref{combinterpBa} 
are in fact equivalent, we need to express $C_\alpha[X;q]$ in terms of
$B_\beta[X;q]$.

%allows us to conclude that 
%$C_\alpha[X;q] = \omega q^{|\alpha| - \ell(\alpha)} B_{\overleftarrow\alpha}[X;1/q]$
%and similarly, 
%$B_\alpha[X;q] = \omega q^{|\alpha| - \ell(\alpha)} C_{\overleftarrow\alpha}[X;1/q]$.

\begin{lem}
\label{www}
Let $\gamma, \alpha$ be compositions with $\gamma \le \alpha$.  Then
\begin{align*}
\sum_{{\beta \atop \gamma \leq \beta \leq \alpha }}
(-1)^{ \ell(\alpha) - \ell(\beta) } q^{ \ell(\alpha) - \ell(\beta)+\doff_{\beta}( DP(\gamma) ) -
\doff_{\overleftarrow{\alpha}}(DP(\overleftarrow{\beta}))} = 
\begin{cases} 0 \text{ if $\gamma < \alpha$}\\
1 \text{ if $\gamma = \alpha$}
\end{cases}.
\end{align*}
\end{lem}
\begin{proof}
First assume that $\gamma < \alpha$
and consider the difference of the
descent sets $\Des(\gamma)\backslash \Des(\alpha) = \{ i_1, i_2, \ldots, i_d\}$.
There are $2^d$ compositions $\beta$ such that $\gamma \leq \beta \leq \alpha$
and we will pair them up with a sign reversing involution.

If $i_1 \in \Des(\beta)$, let $\tilde{\beta}$ be the composition with $\Des(\tilde{\beta})
= \Des(\beta)\backslash\{i_1\}$ (the terms with $i_1 \in \Des(\beta)$ will match with the
terms $i_1 \notin \Des(\tilde{\beta})$).  There is some $r > 1$ such that
$\alpha_1 = \beta_1$, $\alpha_2 = \beta_2, \ldots, \alpha_r > \beta_r$ because $\Des(\beta)$
contains the descent $i_1$ that is not in $\Des(\alpha)$.  Calculating directly
we have that,
\begin{align*}
\ell(\alpha) - \ell(\beta) &= \ell(\alpha) - \ell(\tilde{\beta}) - 1 \\
\doff_\beta(DP(\gamma)) &= \doff_{\tilde{\beta}}(DP(\gamma)) + r \\
-\doff_{\overleftarrow{\alpha}}(DP(\overleftarrow{\beta})) &=
- \doff_{\overleftarrow{\alpha}}(DP(\overleftarrow{\tilde{\beta}}))-r+1.
\end{align*}
Therefore
$$
\ell(\alpha) - \ell(\beta)+\doff_\beta(DP(\gamma))-
\doff_{\overleftarrow{\alpha}}(DP(\overleftarrow{\beta}))=
\ell(\alpha) - \ell(\tilde{\beta})+\doff_{\tilde{\beta}}(DP(\gamma))
- \doff_{\overleftarrow{\alpha}}(DP(\overleftarrow{\tilde{\beta}}))
$$
and the signs of these terms are different in the sum.  This provides
a sign reversing involution and hence those terms with $\touch(D) < \alpha$
sum to $0$ matching those terms with $\Des(\beta)$ which include the smallest
descent.

Now for those terms with $\touch(D) = \alpha$ we have $\gamma= \beta = \alpha$
and
$$\ell(\alpha) - \ell(\beta) + \doff_\beta(DP(\gamma)) - \doff_{\overleftarrow{\alpha}}(DP(
\overleftarrow{\beta})) = 0$$
since $\doff_\alpha(DP(\alpha)) = \pchoose{\ell(\alpha)}{2}$.  
\end{proof}

\begin{thrm}
\label{CatoBaexpansion}
For $n \geq 0$ and for any composition $\alpha \models n$,
\begin{equation} \label{CatoBa}
C_{\alpha}[X;q] = \sum_{\beta \leq \alpha} (-q)^{\ell(\alpha)-\ell(\beta)} 
q^{-\doff_{\overleftarrow\alpha}(DP(\overleftarrow\beta))} B_\beta[X;q]~.
\end{equation}
\end{thrm}
\begin{proof} 
By Theorem \ref{BatoCaexpansion} we have
\begin{align*}
\sum_{\beta \leq \alpha} &(-q)^{ \ell(\alpha)-\ell(\beta) } 
q^{ -\doff_{\overleftarrow\alpha}(DP(\overleftarrow\beta)) } B_\beta[X;q]  \\
&=\sum_{\beta \leq \alpha} (-q)^{ \ell(\alpha)-\ell(\beta) } 
q^{ -\doff_{\overleftarrow\alpha}(DP(\overleftarrow\beta)) } 
\sum _{\gamma \le \beta} q^{ \text{doff}_{\beta} (DP(\gamma)) } C_{\gamma}[X;q] \\
 &= \sum_{\gamma} C_{\gamma}[X;q]
\sum_{\beta \atop  \gamma \le \beta \le \alpha}
(-q)^{\ell(\alpha)-\ell(\beta)} 
q^{-\doff_{\overleftarrow\alpha}(DP(\overleftarrow\beta))} 
q^{\text{doff}_{\beta} (DP(\gamma))}\,,
\end{align*}
and our claim now follows by Lemma~\ref{www}.
\end{proof}

\begin{thrm} \label{combequiv}
Conjecture \ref{combinterpCa} is true if and only if Conjecture \ref{combinterpBa} is true. 
\end{thrm}

\begin{proof}  
Theorem~\ref{BatoCaexpansion} gives $B_\alpha[X;q]$ in terms 
of $C_\beta[X;q]$, to which we apply $\nabla$:
\begin{align*}
\nabla(B_{\alpha}[X;q]) &= 
\sum_{\beta \leq \alpha} q^{\doff_\alpha(DP(\beta))} \nabla(C_\beta[X;q])
\,.
\end{align*}
Given that Conjecture \ref{combinterpCa} holds, we then have that
\begin{align*}
\nabla(B_{\alpha}[X;q]) &= 
\sum_{\beta \leq \alpha} \, 
\sum_{ D \atop \touch(D) = \beta} \, \sum_{w \in \WP_{D}} t^{\area(D)} 
q^{\dinv(w)+\doff_\alpha(DP(\beta))} x^w\\
&= 
\sum_{D \leq DP(\beta)} \sum_{w \in \WP_{D}} t^{\area(D)} q^{\dinv(w)+\doff_\alpha(DP(\beta))} x^w~.
\end{align*}

On the other hand, assuming Conjecture \ref{combinterpBa} holds,
Theorem~\ref{CatoBaexpansion} gives 
\begin{align}
\nabla(C_{\alpha}[X;q]) &= \sum_{\beta \leq \alpha} (-q)^{\ell(\alpha)-\ell(\beta)} 
q^{-\doff_{\overleftarrow\alpha}(DP(\overleftarrow\beta))} \nabla(B_\beta[X;q])\nonumber\\
&= \sum_{\beta \leq \alpha} \, \sum_{D \leq DP(\beta)} \, \sum_{w \in \WP_D}
t^{\area(D)} (-q)^{\ell(\alpha)-\ell(\beta)} q^{\dinv(w)+\doff_\beta(D)-
\doff_{\overleftarrow{\alpha}}(DP(\overleftarrow{\beta}))} x^w~.\label{tocancel}
\end{align}
For a Dyck path $D$ in the sum, let $\gamma = \touch(D)$.  Note that
$D \leq DP(\beta)$ implies that $\gamma \leq \beta \leq \alpha$, and
thus we may rearrange sums as follows:
\begin{align*}
\sum_{\beta \leq \alpha} \, &\sum_{D \leq DP(\beta)}
(-q)^{\ell(\alpha) - \ell(\beta)} q^{\doff_\beta(D) - 
\doff_{\overleftarrow{\alpha}}(DP(\overleftarrow{\beta}))}\\
&= \sum_{\gamma \leq \alpha} \sum_{D \atop \touch(D) = \gamma}
\sum_{\beta \atop \gamma \leq \beta \leq \alpha}
(-1)^{\ell(\alpha) - \ell(\beta)} q^{\ell(\alpha) - \ell(\beta)+\doff_\beta(D) - 
\doff_{\overleftarrow{\alpha}}(DP(\overleftarrow{\beta}))}\,.
\end{align*}
Lemma~\ref{www} allows us to conclude that 
\eqref{tocancel} reduces to Conjecture~\ref{combinterpCa}.
\end{proof}

\begin{cor}
Conjecture \ref{belowCI} is true if and only if Conjecture \ref{touchCI} is true.
\end{cor}

%\begin{section}{A second look at an automorphism on Dyck paths}
%In \cite{Hag08}, an automorphism on the set of Dyck paths of length $n$
%is described which has the property that $\phi: DP_n \rightarrow DP_n$
%such that $\area(\phi(D)) = \dinv(D)$ and $bounce(\phi(D)) = \area(D)$.  We
%revisit this map by looking closely at the recurrence proven in \cite{GXZ10}
%which is used to demonstrate Conjecture \ref{touchCI}.

%In that paper it was proven that
%$$\left< e_n[X], \nabla( C_\alpha[X;q]) \right> =
%q^{\ell-1} t^{\alpha_1-1}
%\left< e_{n-1}[X], \nabla( \Cop_{\alpha_2} \Cop_{\alpha_3} \cdots \Cop_{\alpha_\ell} \Bop_{\alpha_1-1}(1)
%\right>~.$$

%\todo{I am not sure if we want to visit this recurrence in this paper.  It slightly goes
%off the track of the rest of the paper}

%\end{section}

%\bibliography{/home/jhaglund/books/qtcat/Hall}
\bibliographystyle{amsalpha} 
%\bibliography{Hall}

%\end{document}

\end{document}